\documentclass [12pt]{article}
\usepackage {amsfonts}
\usepackage {mathrsfs}
\usepackage {amsmath}
\usepackage {latexsym}
\usepackage {amssymb}
\usepackage {amsthm}
\usepackage {amscd}
\date{}
\title{\bf Nichols algebras over classical Weyl groups(II)
\author{\small  Weicai Wu \\
\small $a$. School of Mathematics,  Hunan Institute of Science and Technology\\
\small  Yueyang 414006,   \ P.R. China \\
\small {\tt Emails: weicaiwu@hnu.edu.cn} }}

\begin{document}
\newtheorem{Proposition}{Proposition}[section]
\newtheorem{Theorem}[Proposition]{Theorem}
\newtheorem{Definition}[Proposition]{Definition}
\newtheorem{Corollary}[Proposition]{Corollary}
\newtheorem{Lemma}[Proposition]{Lemma}
\newtheorem{Example}[Proposition]{Example}
\newtheorem{Remark}[Proposition]{Remark}
\newtheorem{Conjecture}[Proposition]{Conjecture}
\maketitle %\addtocounter{section}{-1}

\begin {abstract} It is shown that except in three cases conjugacy classes of classical Weyl groups $W(B_{n})$ and $W(D_{n})$ are of type ${\rm D}$. This proves that Nichols algebras of irreducible Yetter-Drinfeld modules over the classical Weyl groups $\mathbb W_{n}$ (i.e. $H_{n}\rtimes  \mathbb{S}_{n}$) are infinite dimensional, except the class of type $(2, 3),(1^{2}, 3)$ in $\mathbb S_{5}$, and $(1^{n-2},  2)$ in $\mathbb S_{n}$ for $n >5$.
\vskip.2in
\noindent {\em 2020 Mathematics Subject Classification}: 16W50, 16T05\\
{\em Keywords}: Rack; classical Weyl groups; Nichols algebras.
\end {abstract}

\section {Introduction}\label {s0}
This work contributes to the classification of finite-dimensional pointed Hopf algebras over an algebraically closed field  of characteristic $0$. Those over finite abelian groups are relatively clear due to the work of Andruskiewitsch, Schneider, Heckenberger,see \cite {AS10},\cite {He06a},\cite {He06b}, etc. However not much is known about those over nonabelian groups. Given a group $ G$, an important step to classify all finite-dimensional pointed
Hopf algebras $H$ with group-like $G(H) = G$ is to determine all the pairs $(\mathcal O, \rho )$ such that the associated Nichols algebra $  \mathfrak B(\mathcal O, \rho )$ is finite-dimensional. In general, it is often useful to discard those pairs such that $\dim  \mathfrak B (\mathcal O, \rho ) = \infty$. There are properties of the conjugacy class $\mathcal O$ that imply that $\dim  \mathfrak B (\mathcal O, \rho )  = \infty$ for any $\rho $, one of which is the property of being of type {\rm D}. This is useful since it reduces the computations to operations inside the group and avoids hard calculations of generators and relations of the corresponding Nichols algebra.

Nichols algebras of braided vector spaces  $(\mathbb CX,  c_{q})$,  where $X$ is a
rack and $q$ is a 2-cocycle in $X$, were studied in \cite {AG03}.  It was shown
\cite {AFGV08,  AFZ09,  AZ07} that Nichols algebras  $\mathfrak B({\mathcal O}_
\sigma, \rho)$ over symmetry groups  are infinite dimensional,  except in a
number of remarkable cases corresponding to ${\mathcal O}_{\sigma}$. The present paper aims to classify finite-dimensional Nichols algebra over the classical Weyl groups.
Obviously, $S_{n}$ is isomorphic to
Weyl group of $A_{n-1}$ for $n>1$. We give the following  result for completeness.

{\bf Theorem}  (\cite [Th. 1.1] {AFGV08})
Let $n\ge 5$.  Let $\sigma\in \mathbb S_{n}$ be of type $(1^{m_{1}},2^{m_{2}},\dots,n^{m_{n}})$,
let $\mathcal O_\sigma $ be the conjugacy class of $\sigma$ and let $\rho=(\rho,V) \in
\widehat{{S_{n} } ^{\sigma}}$.
If $\dim \mathfrak B ( \mathcal O_\sigma, \rho) < \infty$, then
the type of $\sigma$ and $\rho$ are in the following list:
\renewcommand{\theenumi}{\roman{enumi}}   \renewcommand{\labelenumi}{(\theenumi)}
\begin{enumerate}
\item $(1^{n_1}, 2)$, $\rho_1 = {\rm sgn}$ or $\epsilon$, $\rho_2 ={\rm sgn}$.
\item $(2, 3)$ in $\mathbb S_{5}$, $\rho_2 ={\rm sgn}$, $\rho_3= \overrightarrow{\chi_{0}}$.
\item $(2^3)$ in $\mathbb S_{6}$, $\rho_2=\overrightarrow{\chi_{1}}\otimes\epsilon$ or
$\overrightarrow{\chi_{1}}\otimes {\rm sgn}$.
\end{enumerate}

In \cite{ZZ12} they showed  that except in several cases Nichols algebras
of irreducible Yetter-Drinfeld ({\rm YD} in short) modules over classical Weyl groups $A
\rtimes \mathbb S_{n}$ supported by $\mathbb S_{n}$ are infinite dimensional. One of the present author \cite{TWZ21} showed  that except in several cases conjugacy classes of classical Weyl groups $W(B_{n})$ and $W(D_{n})$ are of type {\rm D}, they prove that except in several cases Nichols algebras of irreducible Yetter-Drinfeld modules over the classical Weyl groups are infinite dimensional.
However, the classification has not been completed for Nichols algebras over
classic Weyl groups $W(B_{n})$ and $W(D_{n})$.

Note that $\mathbb Z_{2}^{n} \rtimes \mathbb S_{n}$
is isomorphic to Weyl groups $W(B_{n})$ and $W(C_{n})$ of $B_{n}$ and
$C_{n}$ for $n>2$. If $K_{n}= \{a \in
     \mathbb Z_{2}^{n} \mid \ a = ({a_{1}}, {a_{2}}, \cdots, {a_{n}}) $  $\hbox {
with   } a_{1} +a_{2} + \cdots + a_{n} \equiv 0  \ ({\rm mod } \ 2) \}$,
then $K_{n} \rtimes \mathbb S_n$ is isomorphic to Weyl group
$W(D_{n})$ of $D_{n}$ for $n>3$ and $ K_{n} \rtimes  \mathbb{S}_{n}$  is a subgroup of
$\mathbb Z_{2}^{n} \rtimes \mathbb S_{n}$. Let $\mathbb W_{n}$ denote $H_{n}\rtimes  \mathbb{S}_{n}$ throughout this paper, where
$H_{n}:=  K_{n}$  or $\mathbb Z_{2}^{n}$. Here is our main result.

{\bf Theorem 4.3} Assume $n\ge 5$. Let $\sigma\in\mathbb S_{n}$ be of type $(1^{m_{1}},2^{m_{2}},\dots,n^{m_{n}})$ and $a\in H_{n}$ with $a\sigma\in \mathbb W_{n}$ and $\sigma\not=1$. If {\rm dim} $\mathfrak B(\mathcal{O}_{a\sigma}^{\mathbb W_{n}}, \rho)<\infty $,  then the type of $\sigma$ belongs to one in the following list:
\renewcommand{\theenumi}{\roman{enumi}}   \renewcommand{\labelenumi}{(\theenumi)}
\begin{enumerate}
\item $(2, 3)$;
\item $(1^{2}, 3)$;
\item $(1^{n-2},  2)$ for $n >5$ with $a_{i} = a_{j}$ when $\sigma(i) = i $ and $\sigma(j)=j$.
\end{enumerate}

Note that $\mathbb{S}_{n}$ acts on
$H_{n}$ as follows: for any $a\in H_{n}$ with $a = ({a_{1}}, {a_{2}},
\cdots, {a_{n}})$ and $\sigma, \tau \in \mathbb{S}_{n}$,
$\sigma (a) :=( {a_{\sigma^{-1}(1)}}, {a_{\sigma^{-1}(2)}}, \cdots,{a_{\sigma^{-1}(n)}} ).$

It is clear that
 \begin {eqnarray*}
&& (a, \sigma )^{- 1} = (- (a _{\sigma (1)}, a_{\sigma (2)}, \cdots  , a_{\sigma (n)}), \sigma^{ - 1}) = (- \sigma ^{ - 1}(a), \sigma ^{- 1}), \\
&&(b, \tau  )(a, \sigma )(b, \tau  )^{- 1} = (b + \tau  (a) -  \tau  \sigma \tau ^{ - 1}(b), \tau  \sigma \tau ^{ - 1}).
\end  {eqnarray*}
Let $a \sigma$ denote $(a, \sigma)$ in short, sometimes.  When we say  type
{\rm D} we mean a rack or conjugacy class of type {\rm D} and not a Weyl group of type
$D_{n}$. The other notations are the same as in \cite {TWZ21}.

In this paper we prove  that except in three cases conjugacy classes of classical Weyl groups $W(B_{n})$ and $W(D_{n})$ are of type ${\rm D}$. This proves that Nichols algebras of irreducible Yetter-Drinfeld modules over the classical Weyl groups $\mathbb W_{n}$ are infinite dimensional, except the class of type $(2, 3),(1^{2}, 3)$ in $\mathbb S_{5}$, and $(1^{n-2},  2)$ with $a_{i} = a_{j}$ when $\sigma (i) = i $ and $\sigma (j)=j$ in $\mathbb S_{n}$ for $n >5$.

The work is organized as follows. In Section 2 we provide some preliminaries and set our notations.
In Section 3  we examine the rack and
conjugacy classes of $\mathbb W_{n}$.
In  Section 4 we classify Nichols algebras of irreducible {\rm YD} modules over the classical Weyl groups, finally we give a Conjecture: Assume $n\geq5$. then {\rm dim} $\mathfrak B(\mathcal{O}_{a\sigma}^{\mathbb W_{n}}, \rho)=\infty $ for $a\sigma\in \mathbb W_{n}$ with $\sigma\in\mathbb S_{n}, \sigma\not=1$ and $a\in H_{n}$.

\section{ Preliminaries}\label {s1}
Let {\footnotesize $\mathbb{K} $} be an algebraically closed field of characteristic zero and $(V, C)$
 a braided vector space, i.e.   $V$ is
 a  vector space over {\footnotesize $\mathbb{K} $} and $C \in {\rm Aut} (V \otimes V)$ is a solution
 of the braid equation $   ({\rm id} \otimes C_{V,V}) (C_{V,V}\otimes {\rm id})({\rm id} \otimes C_{V,V}) = (C_{V,V}\otimes {\rm id})({\rm id} \otimes C_{V,V})( C_{V,V}\otimes{\rm id})$.
 $C_{i, i+1}:=  {\rm id} ^{\otimes i-1}     \otimes C_{V, V} \otimes {\rm id} ^{\otimes n-i-1}$.
Let $S_{m}\in {\rm End}_{\mathbb{ K} }(V^{\otimes m})$
and $S_{1, j}\in {\rm End}_{\mathbb{ K}}(V^{\otimes (j+1)})$ denote the maps $S_{m}=\prod \limits_{j=1}^{m-1}({\rm id}^{\bigotimes m-j-1}  \bigotimes S_{1, j})$ ,
 $S_{1, j}={\rm id}+C_{12}^{-1}+C_{12}^{-1}C_{23}^{-1}+\cdots+C_{12}^{-1}C_{23}^{-1}\cdots C_{j, j+1}^{-1}$
(in leg notation) for $m\geq 2$ and $j\in \mathbb N$. Then the subspace
$S=\bigoplus \limits _{m=2}^{\infty}  {\rm ker} S_{m}$ of the tensor $T(V)=\bigoplus \limits _{m=0}^{\infty}V^{\otimes m}$
is a two-sided ideal,  and algebra $\mathfrak B(V)=T(V)/S$ is termed the Nichols algebra associated to $(V, C)$.

For  $s\in G$ and  $(\rho,  V) \in  \widehat {G^s}$,  here is a
precise description of the {\rm YD} module $M({\mathcal O}_s,
\rho)$,  introduced in \cite {Gr00}. Let $t_1 = s, t_2,   \cdots,
t_{m}$ be a numeration of ${\mathcal O}_s$,  which is a conjugacy
class containing $s$,   and let $g_i\in G$ such that $g_i \rhd s :=
g_i s g_i^{-1} = t_i$ for all $1\le i \le m$. Then $M({\mathcal
O}_s,  \rho) = \oplus_{1\le i \le m}g_i\otimes V$. Let $g_iv :=
g_i\otimes v \in M({\mathcal O}_s, \rho)$,  $1\le i \le m$,  $v\in V$.
If $v\in V$ and $1\le i \le m $,  then the action of $h\in G$ and the
coaction are given by
\begin {eqnarray} \label {e0.11}
\delta(g_iv) = t_i\otimes g_iv,  \qquad h\cdot (g_iv) =
g_j(\nu _i(h)\cdot v),  \end {eqnarray}
 where $hg_i = g_j\nu _i(h)$,  for
unique  $1\le j \le m$ and $\nu _i(h)\in G^s$. The explicit formula for
the braiding is then given by
\begin{equation} \label{yd-braiding}
C(g_iv\otimes g_jw) = t_i\cdot(g_jw)\otimes g_iv =
g_{j'}(\nu _j (t_i)\cdot w)\otimes g_iv\end{equation} for any $1\le i, j\le
m$,  $v, w\in V$,  where $t_ig_j = g_{j'}\nu _j(t_i)$ for unique $j'$,  $1\le
j' \le m$ and $\nu _j(t_i) \in G^s$. Let $\mathfrak{B} ({\mathcal O}_s,
\rho )$ denote $\mathfrak{B} (M ({\mathcal O}_s,  \rho ))$.
$M({\mathcal O}_s,  \rho )$ is a simple {\rm YD} module (see \cite
 {AZ07}).

We briefly recall the definition and main properties of racks;
see \cite {AG03} for details, more information and bibliographical references.
\textit{A rack} is a pair $(X, \rhd )$,  where $ X$
is a non-empty set and  $\rhd :  X \times X \rightarrow  X$ is an operation such that
  $ x\rhd x = x$, $x\rhd (y \rhd z) = (x\rhd y)\rhd (x\rhd z)$  and $\phi _x$
is invertible for any $x, y, z
\in  X$, where $\phi _x$ is a map from $X$ to $X$ sending $y$ to $x\rhd y$  for any $x, y\in X.$
For example, $({\mathcal O}_s^G, \rhd )$ is a rack  with $x\rhd y:= x y x^{-1}.$

If $R$ and $S$ are two subracks of $X$  with $R\cup S = X$, $R\cap S =\emptyset$,  $x\rhd y \in S$, $y\rhd x \in R$,  for any $x\in R, y\in S,$ then $R\cup S$ is called a \textit{decomposition of subracks} of $X$. Furthermore, if there exist $a\in R$, $b\in S$ such that ${\rm sq} (a, b) := a\rhd (b \rhd (a \rhd b)) \not= b$, then $X$ is called to be of \textit{type  D}. Notice that if a rack $Y$ contains a subrack  $X$ of type ${\rm D}$, then $Y$  is also of  type ${\rm D}$ (see \cite {AFGV08}).

\section{Conjugacy classes of $\mathbb W_{n}$}\label {s2}
In this section we  prove that except in three cases conjugacy classes of classical Weyl groups $\mathbb W_n$ are of type {\rm D}.

\begin {Lemma}\label {1.1} (\cite [Lem. 3.1(iii)] {TWZ21}) Let $G= \mathbb W_n$ and $a\sigma, b\tau\in G$.
If $\tau\sigma=\sigma\tau$,  then ${\rm sq }(a\sigma, b\tau) = b\tau$ if and only if
\begin {eqnarray}\label {e1.1}  a+\sigma\tau(a)+\sigma\tau^{2}(a)+\tau(a)
=b+\sigma(b)+\sigma^{2}\tau(b)+\sigma \tau(b).\end {eqnarray}
furthermore,
If $\tau\sigma=\sigma\tau$ and $\sigma^{2}=\tau^{2}=1$,  then ${\rm sq }(a\sigma, b\tau) = b\tau$ if and only if
\begin {eqnarray}\label {e1.2}  a+\sigma\tau(a)+\sigma(a)+\tau(a)
=b+\sigma(b)+\tau(b)+\sigma \tau(b).\end {eqnarray}
\end {Lemma}

\begin {Lemma}\label{1.2} Let $a=(a_{1},\ldots,a_{n}),b=(b_{1},\ldots,b_{n})$ and
$\sum\limits_{i=1}^{n}a_{i}\equiv\sum\limits_{i=1}^{n}b_{i}(mod\ 2)$. Assume
$a\sigma=(a^{(1)}\sigma_{1})(a^{(2)}\sigma_{2})\cdots(a^{(r)}\sigma_{r})$
and
$ b\tau=(b^{(1)}\tau_{1})(b^{(2)}\tau_{2})\cdots(b^{(r)}\tau_{r})$
are independent sign cycle decomposition (see \cite [Theo. 5.4] {ZZ12}) of $a\sigma, b\tau\in \mathbb W_n$, respectively.
If the lengths of
$\sigma_{i},\tau_{i}$ are the same for all $1\leq i\leq r$,  then $a\sigma$ and $b\tau$ are conjugate in $\mathbb W_{n}$.
\end {Lemma}
\noindent {\it Proof.} Assume $\mathbb W_{n}:=\mathbb Z_{2}^{n} \rtimes \mathbb S_{n}$. By $\sum\limits_{i=1}^{n}a_{i}\equiv\sum\limits_{i=1}^{n}b_{i}(mod\ 2)$,
we have $\sum\limits_{i=1}^{n}(a^{(1)}_{i}+a^{(2)}_{i}+\cdots+a^{(r)}_{i})\equiv
\sum\limits_{i=1}^{n}(b^{(1)}_{i}+b^{(2)}_{i}+\cdots+b^{(r)}_{i})(mod\ 2)$. It is enough to show that
$\sum\limits_{i=1}^{n}a^{(j)}_{i}\equiv
\sum\limits_{i=1}^{n}b^{(j)}_{i}(mod\ 2)$ for $\forall\ 1\leq j\leq r$ after rearranging.
Then we know $a^{(j)}\sigma_{j}$
and $b^{(j)}\tau_{j}$ are conjugate for $\forall\ 1\leq j\leq r$ by \cite [Lemm. 5.3] {ZZ12} since
the lengths of
$\sigma_{i},\tau_{i}$ are the same for all $1\leq i\leq r$. Consequently, $a\sigma$
and $b\tau$ are conjugate by \cite [Theo. 5.4] {ZZ12}.

Assume $\mathbb W_{n}:=K_{n} \rtimes \mathbb S_{n}$, it is clear by \cite [Lemm. 5.6] {ZZ12}.
\hfill $\Box$

\begin {Proposition}\label {1.3} If $\sigma$  is of type $(1,  2^{2})$,  then
$\mathcal O_{a\sigma}^{\mathbb W_{5}}$ is of type ${\rm D}$ for all $a \in H_{5}$.
\end {Proposition}
\noindent {\it Proof.} Let $\sigma = (1\ 2)(3\ 4)$ without lost generality, take
$\tau =(1\ 3)(2\ 4)$.

{\rm (i)} Assume $\mathbb W_{5}:=\mathbb Z_{2}^{5} \rtimes \mathbb S_{5}$.
If $a =(a_{1}, a_{2}, a_{3}, a_{4}, 1) $ and $\sum\limits_{i=1}^{4}a_{i}\equiv0 \mbox { or }1(mod\ 2)$,  let $b =(b_{1}, b_{2}, b_{3}, b_{4}, 0) $ and $\sum\limits_{i=1}^{4}b_{i}\equiv1\mbox { or }0(mod\ 2)$, respectively;
if $a =(a_{1}, a_{2}, a_{3}, a_{4}, 0) $ and $\sum\limits_{i=1}^{4}a_{i}\equiv0\mbox { or }1(mod\ 2)$,  let $b =(b_{1}, b_{2}, b_{3}, b_{4}, 1) $ and $\sum\limits_{i=1}^{4}b_{i}\equiv1\mbox { or }0(mod\ 2)$, respectively. Let $R:=\mathbb Z_{2}^{5} \rtimes\sigma \cap  \mathcal O_{a \sigma}^{\mathbb W _{5}} $   and
$S:=\mathbb Z_{2}^{5} \rtimes \tau \cap\mathcal O _{a \sigma}^{\mathbb W _{5}}$.

{\rm (ii)} Assume $\mathbb W_{5}:=K_{5} \rtimes \mathbb S_{5}$.
If $a =(a_{1}, a_{2}, a_{3}, a_{4}, 1) $ and $\sum\limits_{i=1}^{4}a_{i}\equiv1(mod\ 2)$,
let $b =(b_{1}, b_{2}, b_{3}, b_{4}, 0) $ and $\sum\limits_{i=1}^{4}b_{i}\equiv0(mod\ 2)$;
if $a =(a_{1}, a_{2}, a_{3}, a_{4}, 0) $ and $\sum\limits_{i=1}^{4}a_{i}\equiv0(mod\ 2)$,
let $b =(b_{1}, b_{2}, b_{3}, b_{4}, 1) $ and $\sum\limits_{i=1}^{4}b_{i}\equiv1(mod\ 2)$.
Let $R:=K_{5} \rtimes\sigma \cap  \mathcal O_{a \sigma}^{\mathbb W _{5}} $   and
$S:=K_{5} \rtimes \tau \cap\mathcal O _{a \sigma}^{\mathbb W _{5}}$.

Thus $a \sigma $ and $b \tau$ are conjugate by lemma \ref {1.2}. Obviously,  $\sum\limits_{i=1}^{4}a_{i}\neq \sum\limits_{i=1}^{4}b_{i}$ implies ${\rm sq}(a\sigma, b\tau)\not=b\tau$ by
\cite [Exam.4.1]{TWZ21}.  It is clear that
$R \cup S$ is a subrack decomposition of   $\mathcal O_{a\sigma}^{\mathbb W_{5}}$, hence it is of type
${\rm D}$. \hfill $\Box$

\begin {Proposition}\label {1.4} If $\sigma$  is of type $(1^{2},  2^{2})$,  then
$\mathcal O_{a\sigma}^{\mathbb W_{6}}$ is of type ${\rm D}$ for all $a \in H_{6}$.
\end {Proposition}
\noindent {\it Proof.} Let $\sigma = (1\ 2)(3\ 4)$ without lost generality, take
$\tau =(1\ 3)(2\ 4)$.
Asuume $a =(a_{1}, a_{2}, a_{3}, a_{4}, a_{5}, a_{6}) ,b =(b_{1}, b_{2}, b_{3}, b_{4}, b_{5}, b_{6}) $.

{\rm (i)} Assume $\mathbb W_{6}:=\mathbb Z_{2}^{6} \rtimes \mathbb S_{6}$.
If $\sum\limits_{i=1}^{4}a_{i}\equiv0(mod\ 2),a_{5}+a_{6}\equiv0\mbox { or }1(mod\ 2)$,
let $\sum\limits_{i=1}^{4}b_{i}\equiv1(mod\ 2),b_{5}+b_{6}\equiv1\mbox { or }0(mod\ 2)$, respectively;
if $\sum\limits_{i=1}^{4}a_{i}\equiv1(mod\ 2),a_{5}+a_{6}\equiv0\mbox { or }1(mod\ 2)$,
let $\sum\limits_{i=1}^{4}b_{i}\equiv0(mod\ 2),b_{5}+b_{6}\equiv1\mbox { or }0(mod\ 2)$, respectively.
Let $R:=\mathbb Z_{2}^{6} \rtimes\sigma \cap  \mathcal O_{a \sigma}^{\mathbb W _{6}} $   and
$S:=\mathbb Z_{2}^{6} \rtimes \tau \cap\mathcal O _{a \sigma}^{\mathbb W _{6}}$.

{\rm (ii)} Assume $\mathbb W_{6}:=K_{6} \rtimes \mathbb S_{6}$.
If $\sum\limits_{i=1}^{4}a_{i}\equiv0(mod\ 2),a_{5}+a_{6}\equiv0(mod\ 2)$,
let $\sum\limits_{i=1}^{4}b_{i}\equiv1(mod\ 2),b_{5}+b_{6}\equiv1(mod\ 2)$;
if $\sum\limits_{i=1}^{4}a_{i}\equiv1(mod\ 2),a_{5}+a_{6}\equiv1(mod\ 2)$,
let $\sum\limits_{i=1}^{4}b_{i}\equiv0(mod\ 2),b_{5}+b_{6}\equiv0(mod\ 2)$.
Let $R:=K_{6} \rtimes\sigma \cap  \mathcal O_{a \sigma}^{\mathbb W _{6}} $   and
$S:=K_{6} \rtimes \tau \cap\mathcal O _{a \sigma}^{\mathbb W _{6}}$.

Thus $a \sigma $ and $b \tau$ are conjugate by lemma \ref {1.2}. Obviously,  $\sum\limits_{i=1}^{4}a_{i}\neq \sum\limits_{i=1}^{4}b_{i}$ implies ${\rm sq}(a\sigma, b\tau)\not=b\tau$ by
\cite [Exam.4.1]{TWZ21}. It is clear that
$R \cup S$ is a subrack decomposition of   $\mathcal O_{a\sigma}^{\mathbb W_{6}}$, hence it is of type
${\rm D}$. \hfill $\Box$

\begin {Proposition}\label {1.5}  {\rm (i)} If $\sigma$  is of type $(2^{3})$,  then
$\mathcal O_{a\sigma}^{\mathbb W_{6}}$ is of type ${\rm D}$ for all $a \in H_{6}$.

 {\rm (ii)} If $\sigma$  is of type $(2^{4})$,  then
$\mathcal O_{a\sigma}^{\mathbb W_{8}}$ is of type ${\rm D}$ for all $a \in H_{8}$.
\end {Proposition}
\noindent {\it Proof.} Similarly to lemma \ref {1.4}. {\rm (i)} Let $\sigma = (1\ 2)(3\ 4)(5\ 6)$, take
$\tau =(1\ 3)(2\ 4)(5\ 6)$.
{\rm (ii)} Let $\sigma = (1\ 2)(3\ 4)(5\ 6)(7\ 8)$, take
$\tau =(1\ 3)(2\ 4)(5\ 6)(7\ 8)$. \hfill $\Box$

\begin {Proposition}\label {1.6} If $\sigma$  is of type $(1^{n-3},  3)$ for $n >5$,  then $\mathcal O_{a\sigma} ^{\mathbb W_{n}}$ is of type ${\rm D}$ for all $a \in H_{n}$ with $n >5$.
\end {Proposition}
\noindent {\it Proof.} Let $\sigma = (1\ 2\ 3)$, $\tau= (4\ 5\ 6)$. We have that
the $1$-th,  $2$-th ,$3$-th, $4$-th,  $5$-th and $6$-th components of (\ref {e1.1}) are
\begin {eqnarray}\label {e1.3}  (0,0,0,a_{4}+a_{6},  a_{4}+a_{5},  a_{5}+a_{6}) =  (b_{1}+b_{3},  b_{1}+b_{2},  b_{2}+b_{3},0,0,0). \end {eqnarray}

{\rm (1)} Assume $\mathbb W_{n}:=\mathbb Z_{2}^{n} \rtimes \mathbb S_{n}$.

{\rm (i)} If $4$-th,  $5$-th, $6$-th components of $a$ are $(0,0,0)$ and $\sum\limits_{i\neq4,5,6}a_{i}\equiv0\mbox { or }1(mod\ 2)$,
let $1$-th,  $2$-th, $3$-th components of $b$ are $(1,0,0)$ and $\sum\limits_{i=4}^{n}b_{i}\equiv1\mbox { or }0(mod\ 2)$, respectively,  then
(\ref {e1.3}) does not hold.

{\rm (ii)} If $4$-th,  $5$-th, $6$-th components of $a$ are $(1,0,0)$ and $\sum\limits_{i\neq4,5,6}a_{i}\equiv0\mbox { or }1(mod\ 2)$,
let $1$-th,  $2$-th, $3$-th components of $b$ are $(0,0,0)$ and $\sum\limits_{i=4}^{n}b_{i}\equiv1\mbox { or }0(mod\ 2)$, respectively,  then
(\ref {e1.3}) does not hold.

{\rm (iii)} If $4$-th,  $5$-th, $6$-th components of $a$ are $(1,1,0)$ and $\sum\limits_{i\neq4,5,6}a_{i}\equiv0\mbox { or }1(mod\ 2)$,
let $1$-th,  $2$-th, $3$-th components of $b$ are $(1,0,0)$ and $\sum\limits_{i=4}^{n}b_{i}\equiv1\mbox { or }0(mod\ 2)$, respectively,  then
(\ref {e1.3}) does not hold.

{\rm (iv)} If $4$-th,  $5$-th, $6$-th components of $a$ are $(1,1,1)$ and $\sum\limits_{i\neq4,5,6}a_{i}\equiv0\mbox { or }1(mod\ 2)$,
let $1$-th,  $2$-th, $3$-th components of $b$ are $(1,1,0)$ and $\sum\limits_{i=4}^{n}b_{i}\equiv1\mbox { or }0(mod\ 2)$, respectively,  then (\ref {e1.3}) does not hold.
Let $R:=\mathbb Z_{2}^{n} \rtimes\sigma \cap  \mathcal O_{a \sigma}^{\mathbb W _{n}} $   and
$S:=\mathbb Z_{2}^{n} \rtimes \tau \cap\mathcal O _{a \sigma}^{\mathbb W _{n}}$.

{\rm (2)} Assume $\mathbb W_{n}:=K_{n} \rtimes \mathbb S_{n}$.

{\rm (i)} If $4$-th,  $5$-th, $6$-th components of $a$ are $(0,0,0)$ and $\sum\limits_{i\neq4,5,6}a_{i}\equiv0(mod\ 2)$,
let $1$-th,  $2$-th, $3$-th components of $b$ are $(1,0,0)$ and $\sum\limits_{i=4}^{n}b_{i}\equiv1(mod\ 2)$,  then (\ref {e1.3}) does not hold.

{\rm (ii)} If $4$-th,  $5$-th, $6$-th components of $a$ are $(1,0,0)$ and $\sum\limits_{i\neq4,5,6}a_{i}\equiv1(mod\ 2)$,
let $1$-th,  $2$-th, $3$-th components of $b$ are $(0,0,0)$ and $\sum\limits_{i=4}^{n}b_{i}\equiv0(mod\ 2)$,  then
(\ref {e1.3}) does not hold.

{\rm (iii)} If $4$-th,  $5$-th, $6$-th components of $a$ are $(1,1,0)$ and $\sum\limits_{i\neq4,5,6}a_{i}\equiv0(mod\ 2)$,
let $1$-th,  $2$-th, $3$-th components of $b$ are $(1,0,0)$ and $\sum\limits_{i=4}^{n}b_{i}\equiv1(mod\ 2)$,  then
(\ref {e1.3}) does not hold.

{\rm (iv)} If $4$-th,  $5$-th, $6$-th components of $a$ are $(1,1,1)$ and $\sum\limits_{i\neq4,5,6}a_{i}\equiv1(mod\ 2)$,
let $1$-th,  $2$-th, $3$-th components of $b$ are $(1,1,0)$ and $\sum\limits_{i=4}^{n}b_{i}\equiv0(mod\ 2)$,  then (\ref {e1.3}) does not hold.
Let $R:=K_{n} \rtimes\sigma \cap  \mathcal O_{a \sigma}^{\mathbb W _{n}} $   and
$S:=K_{n} \rtimes \tau \cap\mathcal O _{a \sigma}^{\mathbb W _{n}}$.

Thus $a \sigma $ and $b \tau$ are conjugate by \cite [Lemm. 5.3] {ZZ12}. Obviously,  (\ref {e1.3}) does not hold implies ${\rm sq}(a\sigma, b\tau)\not=b\tau$. It is clear that
$R \cup S$ is a subrack decomposition of   $\mathcal O_{a\sigma}^{\mathbb W_{n}}$, hence it is of type
${\rm D}$.  \hfill $\Box$

\begin{Theorem}\label{1.8}
Let $n>4$. Let $\sigma\in\mathbb S_n$ be of type $(1^{m_{1}},2^{m_{2}},\dots,n^{m_{n}})$ and $a\in H_{n}$ with $a\sigma\in W_{n}$ and $\sigma\not=1$. If $\mathcal{O}_{a\sigma}^{ \mathbb W_{n}}$ is not of
type {\rm D},  then the type of $\sigma$ belongs to one in the following list:
\renewcommand{\theenumi}{\roman{enumi}}   \renewcommand{\labelenumi}{(\theenumi)}
\begin{enumerate}
\item $(2, 3)$;
\item $(1^{2}, 3)$;
\item $(1^{n-2},  2)$ for $n >5$ with $a_{i} = a_{j}$ when $\sigma (i) = i $ and $\sigma (j)=j$.
\end{enumerate}
\end{Theorem}
\noindent {\it proof.} It follows from Proposition \ref {1.3}, Proposition \ref {1.4}, Proposition \ref {1.5} , Proposition \ref {1.6} and \cite [Theorem 4.1] {TWZ21}. \hfill $\Box$

\section {Nichols algebras  of irreducible {\rm YD} module over $\mathbb W_{n}$}\label {s3}
In this section we show that  except in three  cases Nichols algebras of irreducible {\rm YD} modules over classical Weyl groups $\mathbb W_n$ are infinite dimensional.

We shall use the following facts:

\begin{Lemma} \label {2.1}(\cite [Cor. 8.4] {HS08}) Let $n\in \mathbb N $, $n\ge 3$, and assume $U$ be a {\rm YD} module  over symmetric group $\mathbb{S}_{n}$. If $\mathfrak B (U)$ is  finite dimensional, then $U$ is an irreducible {\rm YD} module over $\mathbb{S}_{n}$.
\end {Lemma}

\begin {Proposition}\label {2.3} (See \cite [Theorem 3.6]{AFGV08} or \cite [Theorem 5.3]{TWZ21})
If $G$ is a finite group and  $\mathcal{O}_{\sigma}^G$ is of type {\rm D},  then {\rm dim} $\mathfrak B(\mathcal{O}_{\sigma}^G,  \rho) = \infty $ for any $\rho \in \widehat { G^\sigma}$.
\end {Proposition}

\begin{Theorem}\label{2.6} Assume $n\ge 5$. Let $\sigma\in\mathbb S_{n}$ be of type $(1^{m_{1}},2^{m_{2}},\dots,n^{m_{n}})$ and $a\in H_{n}$ with $a\sigma\in \mathbb W_{n}$ and $\sigma\not=1$. If {\rm dim} $\mathfrak B(\mathcal{O}_{a\sigma}^{\mathbb W_{n}}, \rho)<\infty $,  then the type of $\sigma$ belongs to one in the following list:
\renewcommand{\theenumi}{\roman{enumi}}   \renewcommand{\labelenumi}{(\theenumi)}
\begin{enumerate}
\item $(2, 3)$;
\item $(1^{2}, 3)$;
\item $(1^{n-2},  2)$ for $n >5$ with $a_{i} = a_{j}$ when $\sigma(i) = i $ and $\sigma(j)=j$.
\end{enumerate}
\end{Theorem}
\noindent {\it proof.} It follows from Theorem \ref {1.8} and Proposition \ref {2.3}. \hfill $\Box$

Moreover, we give the following conjecture:
\begin{Conjecture}\label{2.7} Assume $n\geq5$. Let $a\sigma\in \mathbb W_{n}$ with $\sigma\in\mathbb S_{n}, \sigma\not=1$ and $a\in H_{n}$, then {\rm dim} $\mathfrak B(\mathcal{O}_{a\sigma}^{\mathbb W_{n}}, \rho)=\infty $.
\end{Conjecture}

\section*{Acknowledgment}
I was supported by Hunan Provincial Natural Science Foundation of China (Grant No: 2020JJ5210) and A Project Supported by Scientific Research Fund of Hunan Provincial Education Department (Grant No: 19C0860).

\end {document}